\documentclass[12pt]{article}
\usepackage[english]{babel}
\usepackage{amssymb,amsmath,amsthm}
\usepackage{graphicx}

\usepackage{amsmath,euscript}
\usepackage{amssymb}
\usepackage{amsthm}
\usepackage{enumerate}
\usepackage{amsfonts}
\usepackage{comment}
\usepackage{colortbl}
\usepackage{amssymb,amsmath,array}

\usepackage[enableskew]{youngtab}
\usepackage{latexsym}

\usepackage{amscd}

\usepackage[latin1]{inputenc}
\usepackage{pstricks,pst-node,pst-text,pst-3d}
\usepackage{amsmath}
\usepackage{pgf}

\newtheorem{theorem}{Theorem}[section]
\newtheorem{lemma}[theorem]{Lemma}
\newtheorem{proposition}[theorem]{Proposition}
\newtheorem{corollary}[theorem]{Corollary}

{\theoremstyle{definition}}
{\theoremstyle{definition}}
{\theoremstyle{definition}\newtheorem{definition}[theorem]{Definition}}
{\theoremstyle{definition}}

\begin{document}

\def\C{{\mathbb C}}
\def\T{{\mathbb T}}
\def\N{{\mathbb N}}
\def\Z{{\mathbb Z}}
\def\R{{\mathbb R}}
\def\K{{\mathbb K}}
\def\CC{{\cal C}}
\def\H{{\cal H}}
\def\F{{\cal F}}
\def\X{{\cal X}}
\def\Y{{\cal Y}}
\def\epsilon{\varepsilon}
\def\kappa{\varkappa}
\def\phi{\varphi}
\def\leq{\leqslant}
\def\geq{\geqslant}
\def\re{\text{\tt Re}\,}
\def\ilim{\mathop{\hbox{$\underline{\hbox{\rm lim}}$}}\limits}
\def\dim{\hbox{\tt dim}\,}
\def\ker{\hbox{\tt ker}\,}
\def\supp{\hbox{\tt supp}\,}
\def\Re{\hbox{\tt Re}\,}
\def\ssub#1#2{#1_{{}_{{\scriptstyle #2}}}}

\title{K-theory of locally finite graph $C^*$-algebras}

\author{  Natalia~Iyudu
 }

\maketitle

{\bf Abstract}

{\small

\bigskip

We calculate the K-theory of the Cuntz-Krieger algebra ${\cal O}_E$
associated with an infinite, locally finite graph, via  the
Bass-Hashimoto operator. The formulae we get express the
Grothendieck group and  the Whitehead group in purely graph
theoretic terms.

We consider the category of finite
(black-and-white, bi-directed) subgraphs with certain graph
homomorphisms and construct a continuous functor to abelian groups.
In this category $K_0$ is an inductive limit of $K$-groups of
finite graphs, which were calculated in \cite{MM}.

In the case of an infinite graph with the finite Betti
number we obtain the formula for the Grothendieck group $K_0({\cal O}_E)= {\mathbb
Z}^{\beta(E)+\gamma(E)},\,$ where $\beta(E)$ is the first Betti
number and $\gamma(E)$ is the valency number of the graph $E$. We
note, that in the infinite case the torsion part of $K_0$, which
is present in the case of a finite graph, vanishes. The
 Whitehead group depends only on the first Betti
number: $K_1({\cal O}_E)= {\mathbb Z}^{\beta(E)}$. These allow us to
provide a counterexample to the fact, which holds for finite
graphs, that  $K_1({\cal O}_E)$ is
the torsion free part of $K_0({\cal O}_E)$.

\bigskip

%graph C*-algebra, Grothendieck and Whitehead groups,  Betti number, Bass-Hashimoto operator, %Cuntz-Krieger algebra

}

{\centerline\hrulefill}

\bigskip

\normalsize

\section{Introduction}

We start with defining, how the $C^*$-algebra is associated with a graph, in our setting.
Namely, we explain that we deal with  a  $C^*$-algebras associated with a graph via the Bass-Hashimoto operator.

Let us consider first a non-directed graph $\hat E$, which is
allowed to be infinite, to have loops, multiple edges and sinks.

We assume that the graph is {\it locally finite}, i.e. every vertex
is connected only to finitely many vertices by edges, and that
suppose that $\hat E$ is a connected non-directed graph (=geometrically
connected graph).

We deal with the Cuntz-Krieger algebra ${\cal O}_E$ of the
Bass-Hashimoto operator associated with the graph $E$ (more precisely,
its infinite, locally finite analogue). This operator (operator
$\Phi_E$ defined below) was considered by Hashimoto \cite{H} and
Bass \cite{B} and later studied in \cite{MM}. The algebra known as a
{\it boundary operator algebra} (for example, cf. \cite{Ro, Ro1,
Ro2}) is Morita equivalent to the corresponding algebra associated
with the Bass-Hashimoto operator. Namely, ${\cal O}_E \sim C^*(\delta
{\mathbb E})/{\Gamma}$, where $ {\mathbb E}$ is a universal covering
tree of the graph $E$ and ${\Gamma}$ is a free group of rank
$\beta$, where $\beta$ is the first Betti number of $E$.

 Note, that
the Cuntz-Krieger algebra of the Bass-Hashimoto operator associated with a  graph $E$ (denoted here also by ${\cal O}_E$) should not be
mixed with the Cuntz-Krieger algebra of the operator defined by the
incidence matrix of the graph, as it is done, for example, in
\cite{A,Wa}. These are two different ways to associate a
Cuntz-Krieger algebra with a graph, via different operators. Although
the operators appear very similar, the behaviour of the algebra changes
dramatically.
For example, as it can be seen from \cite{A}, $K_0$ of algebras
associated with finite graphs via the incidence matrix  are far from being defined
by the first Betti number, as it is the case for  algebras
associated with a finite graph via the Bass-Hashimoto operator
\cite{MM}.

Our goal here is to calculate the K-theory of the Cuntz-Krieger
algebra associated with the Bass-Hashimoto operator of an infinite,
locally finite graph, purely in graph theoretic terms, as it was
done in \cite{MM} in the case of a finite graph.

The idea of our calculations is essentially the same as in \cite{MM}: we use the fact that groups $K_0$ and $K_1$ do not change under the finite number of edge contractions. In case of the infinite graph, by the finite number of edge contractions we could arrive to a simpler, but still infinite graph. Namely, we arrive at a rose with outgoing trees. In  the case of $K_0$, considering the presentation of this group as a quotient of an infinitely generated free group, we are able to write down relations for finite subsets of generators, and conclude that not only the Betti number, but also another characteristic of the graph,
which we call the 'valency number',  appears in the formula for $K_0$. In the case of $K_1$,  calculating in the infinitely generated group, we prove, that an element from the kernel of the operator, associated with a graph (whose kernel is $K_1$), should contain only linear combinations of generators, corresponding to the edges of the rose, but not of the outgoing trees. Due to this, the formula for $K_1$, as in the finite case, contains only the Betti number of the graph.

For any infinite, locally finite graph, which is connected, we can
define the first Betti number (cyclomatic number), extending the
usual definition for a finite graph.

\begin{definition}\label{d1}\rm If $\widehat E $ is a finite
geometrically connected graph, then the first Betti number of $\widehat E$ is
$$
\beta(\widehat E)=d_1-d_0+1,
$$
where $d_0$ is the cardinality of the set of vertices and $d_1$
is the cardinality of the set of (geometric) edges.
\end{definition}

Note that for a finite graph with $m$ connected components it would
be
$$
\beta(\widehat E)=d_1-d_0+m.
$$
This number determines the number of cycles in $\widehat E$.

\begin{definition}\label{d2}\rm If $\widehat E$ is an infinite, locally finite
geometrically connected graph, we define the first Betti number
of $\widehat E$ as the limit of the sequence of Betti numbers of finite
subgraphs ${\widehat E}_k$ of $\widehat E$ obtained in the following way: $\widehat{E}_0$ is an
arbitrary connected finite subgraph of $\widehat E$ and for any $n$,
${\widehat E}_{n+1}$ is obtained from ${\widehat E}_n$ by adding to $\widehat{E}_n$ all edges of
$\widehat E$ connected to the vertices of $\widehat{E}_n$ (together with the vertices on the other end of these edges). It will be a finite graph, since  $\widehat E$ is locally finite.

\end{definition}

{\bf Remark 1.} This definition does not depend on the choice of the
subgraph $\widehat{E}_0$, from which the sequence $\{\widehat{E}_n\}$ starts. Indeed, suppose one
starts from another graph $\widehat{E}_0'$.  At some step $n$, one will have
the graph $\widehat{E}_0$ as a subgraph of $\widehat{E}_n'$ and $\widehat{E'}_0$ as a subgraph of
$\widehat{E}_n$ (all vertices and edges will be 'eaten' due to the connectedness
of $\widehat E$). It follows that the sequences $\beta (\widehat{E}_n)$ and
$\beta(\widehat{E'}_n)$ have the same limit: either stabilize on the same
positive integer or both grow to infinity.

{\bf Remark 2.} It is clear that an infinite, locally finite graph with the finite first Betti number, has the shape of a finite graph (with the same first Betti number), with finite or infinite number of  outgoing trees.

\def\ee{\overline{e}}

Now associate for convenience with any graph  $\widehat E$ as above,
an oriented {\it bi-directed} graph
 $E=(E^0,E^1,s,r)$  with the set of
vertices $E^0$, set of edges $E^1$ and maps $s,r$ from $E^1$ to
$E^0$ which determine the source and the range of an arrow
respectively ({\it source} and {\it range} maps). The graph $E$ is obtained from $ \widehat E$
by doubling the edges of $\widehat E$, so that each non-oriented edge of $\widehat E$ gives rise to the pair of edges of $E$, $e$ and $\ee$, equipped with opposite orientations.

For any such finite graph $E$ one can associate a Cuntz--Krieger
$C^*$-algebra ${\cal O}_E$, in the way it is done in \cite{MM}. Namely, there
it is considered a Cuntz--Krieger $C^*$-algebra (as it is defined in the Cuntz--Krieger paper \cite{CK})  associated with the
matrix $A_E$. The
matrix $A_E$ is obtained from the graph as a matrix
of the following  operator (homomorphism of countable direct sum ${\mathbb Z}^{(E^1)}$ of copies of $\mathbb Z$), written as follows in the basis labelled by the  set $E^1$ of edges of $E$:

$$
\Phi_E: {\mathbb Z}^{(E^1)} \to {\mathbb Z}^{(E^1)}: \quad e \mapsto
-\ee+\sum_{e':r(e)=s(e')}e'
\eqno{(*)}
$$

This operator was considered in \cite{H} and \cite{B} in
connection with the study of the Ihara zeta function of a graph, and is called the {\it Bass-Hashimoto operator}.

The entries of the matrix $A_E$ are in $\{0,1\}$. If the graph $E$
is finite, then $A_E$ is an $2n\times 2n$ matrix, where $n$ is the
number of geometric edges of the graph $E$ (=the number of edges of
the corresponding  non-oriented graph $\widehat E$).

Denote by
${\cal O}_E$ the Cuntz--Krieger $C^*$-algebra associated with this matrix
$A_E$ as in \cite{CK}. That is, ${\cal O}_E$ is the $C^*$-algebra generated
by $2n$ partial isometries $\{S_j\}_{j=1}^{2n}$ which act on a
Hilbert space in such a way that their support projections
$Q_i=S_i^*S_i$ and their range projections $P_i=S_iS^*_i$ satisfy
the relations
$$
P_iP_j=0\ \ \text{if $i\neq j$}\ \ \text{and}\ \ Q_i=\sum_{j=1}^{2n}
(A_E)_{ij}P_j\ \ \text{for}\ \ 1\leq j\leq 2n.
$$
This definition of ${\cal O}_E$ surely makes sense for an arbitrary locally finite matrix $A_E$, since the relations still contain  finite sums.

So, we have the following.

\begin{definition}\label{d3}\rm For an infinite, row finite graph
$E=(E^0,E^1,s,r)$,  its $C^*$-algebra
${\cal O}_E$ is generated by partial isometries $\{S_i:i\in E^1 \}$
subject to the relations
$$
S^*_iS_i=\sum_{j\in E^1} (A_E)_{ij}S_jS_j^*.
$$

\end{definition}

Note that there could be another way to associate a $C^*$-algebra with a graph, for example, the one which is considered in \cite{A}. It should be distinguished from that described above. In \cite{A} a $C^*$-algebra of the graph is defined as a Cuntz-Krieger algebra of another operator associated with a graph: the operator is represented by the edges adjacency matrix of an oriented graph.

In the paper \cite{MM} a formula was obtained for $K_0({\cal O}_E)$
depending only on the first Betti number  $\beta(E)$ of the graph $E$, for a finite
graph, namely $K_0({\cal O}_E)={\mathbb Z}^{\beta(E)}\oplus {\mathbb Z} / (\beta(E)-1)\mathbb Z$.
%, where  $\beta(E)$ is a first Betti number of  $E$.
As a consequence, $K_1({\cal O}_E)={\mathbb Z}^{\beta(E)}$, since it is well known in the finite graph setting (\cite{Ror}, \cite{MM}), that $K_1$ is a torsion free part of $K_0$.

It was mentioned in \cite{MM}  that it would be interesting to extend these results to infinite, locally finite graphs.
It turned out that the infinite case does indeed reveal some new phenomenon, which is not present in the finite case.
We do answer this question from \cite{MM}  here
when the  first Betti number $\beta(E)$ of an infinite graph is
finite, and show that in the infinite case $K_0$ does not have a
torsion. Moreover, in the infinite case the formula for $K_0$
involves not only the first Betti number, but also another
combinatorial characteristic of the graph: the {\it valency
number}.

Let us give a precise definition.
Suppose we fix a finite subgraph in $\Gamma$ with the Betti number $\gamma(\Gamma)$.
   (The definition obviously does not depend on this choice.)

  The {\it root vertex} is a root of one of the infinite trees. Its {\it valency}, is the number of outgoing edges, continuing to infinity.

  The {\it branching vertex} is the vertex of a tree with one incoming and more than one outgoing edges, which leads to an infinite path. If $n$ is the number of such outgoing edges, the {\it valency} of that vertex is $n-1$.

  \begin{definition}\label{d4}\rm
  {\it Valency number} $\gamma(\Gamma)$  of an infinite, locally
finite graph $E$ with the finite Betti number $\beta(E)$, is the sum of valencies in all branching vertices and root vertices.

\end{definition}

Note, that in case, the valency number is finite, it coincides with
the number of infinite chains outgoing from the finite subgraph of
$E$ (=the number of infinite ends).
However, in general it is not true. This is demonstrated by the following example.

{\bf Example 1.} Consider the full binary tree $\cal B$.  The valency number $\gamma({\cal B})$ is countable:
all its vertices are branching vertices of valency one, so the valency number is equal to the number of vertices of that tree. The infinite ends however can be enumerated by all sequences of $0$s and $1$s.
If for each vertex, we label the edge going to the right by $0$, and the edge going to the left by $1$, the infinite paths  will be marked by all $0-1$ sequences, so there is a continuum of them.

\begin{theorem}\label{t1} Let $E$ be an infinite, locally finite connected graph
with the finite first Betti number $\beta(E)$ and the valency
number  $\gamma(E)$.
Then
$K_0({\cal O}_E)={\mathbb Z}^{\beta(E) + \gamma(E) }$.
\end{theorem}

Let us give here simple example to demonstrate this theorem.

{\bf Example 2.}
Consider the graph with one loop and one outgoing infinite  chain.
Denote generators of the group $K_0$ associated with a loop by $u$ and
$\bar u$, and those associated with the edges from the chain by $x_1,
\bar x_1, x_2, \bar x_2,... $. The group $K_0$, being the co-kernel of
the operator $T = Id - \Phi_E$ is an abelian group with generators
$u, \bar u, x_1, \bar x_1, x_2, \bar x_2,... $ and relations $Tu,
T\bar u, Tx_1, T\bar x_1, Tx_2, T\bar x_2,... $. By  definition
of the operator $T$, $\,Tu=u-(u+x_1), T\bar u=\bar u -(\bar u+x_1),
Tx_1=x_1-x_2, Tx_2=x_2-x_3,..., T\bar x_1=\bar x_1-(u+\bar u), T\bar
x_2=\bar x_2-\bar x_1, T\bar x_3=\bar x_3-\bar x_2...$. The group $coker
T$ is obviously isomorphic to the group generated by $u, \bar u,
x_1, \bar x_1$, subject to relations $x_1=0, \bar x_1=u+\bar u$,
which means that $K_0({\cal O}_E)={\mathbb Z} \oplus {\mathbb Z} =
{\mathbb Z}^{\beta+\gamma}$.

To emphasize    the nature of the new infinite phenomenon of torsion vanishing for $K_0$, we
would like to give one more example here.

{\bf Example 3.}
Let us have two graphs:
finite graph $G_1 - $  n-rose with finite outgoing path, consisting
of edges $x_1, x_2, x_3$, and infinite graph  $G_2 -$  n-rose with an
infinite outgoing chain, consisting of edges $y_1, y_2, y_3,...$. In
$K_0({\cal O}_{G_1})$ the relation corresponding to the last edge
$x_3$ is $x_3=0$, for others in the path, these are $x_2=x_3,
x_1=x_2$. So, in the quotient of the group we have $x_1=x_2=x_3=0$,
which gives the torsion part of the group.

In the case of the infinite path in the graph $G_2$, where there is no
'last' edge, the relations are $x_1=x_2=x_3=...$, so we get just one
extra variable out of any infinite path, and this variable is non
zero.

Finally, in section 4, we calculate $K_1({\cal O}_E)$ and express it
in terms of the first Betti number.

\begin{theorem}\label{t2}
 Let $E$ be an infinite, locally finite connected graph
with the finite first Betti number $\beta(E)$, and ${\cal O}_E$ is the associated
(via the Bass-Hashimoto operator) $C^*$- algebra.
Then $K_1({\cal O}_E)={\mathbb Z}^{\beta(E)}$.

\end{theorem}

As a consequence of our results for $K_0$ and $K_1$
it turns out that $K_1({\cal O}_E)$ is no longer a torsion
free part of $K_0({\cal O}_E)$ in the infinite graph setting, as it
is the case for finite graphs, as shown in \cite{Ror} or \cite{MM}.
So, we have proved the following corllary.

\begin{corollary}\label{cor1}\rm
 The $K_0$ group of any locally finite infinite graph $E$ is the direct limit of groups corresponding to finite subgraphs from the  category $\cal E$.
\end{corollary}

\section{Category of  black-and-white bi-directed graphs and a functor to abelian groups}

Let $E$ be an infinite, locally finite bi-directed graph as above. We
define here a category $\cal E$  of black-and-white 'subgraphs' of $E$.

Let us note that this section is not necessary for the proof of the main results, so it can be considered as an independent part of the paper aiming to make a link and comparison to the category, considered by P.Ara \cite{A} in the setting of another type of $C^*$-algebras associated with a graph (via a transcendency matrix)
and with other graph categories.
For example, the objects of our category, which we call 'black-and-white graphs', are the same as 'graphs with flags' used in the paper by Yu.Manin and D.Borisov \cite{MB} (flags there are our white edges), but there the set of morphisms is different.
Analogously to \cite{A}, we prove the basic property of the functor from our category to abelian groups.
This section also provides additional insight into  the nature of our calculations for the main results.

To construct  an object of $\cal E$, choose an arbitrary set
$\Omega$ of vertices of $E$. The corresponding graph will contain all
edges of $E$ starting and ending on vertices from $\Omega$. If the
edge starts (ends) on a vertex from $\Omega$, but ends (starts)
outside, it will be called white, otherwise black. Thus the objects of
$\cal E$ are certain (black-and-white) subgraphs of $E$. In
particular,  if we choose the whole set of vertices of $E$, as
$\Omega$, we arrive at the whole graph $E$ (with only black edges
present), as an element of $\cal E$.

 Now we define the set of graph
homomorphisms between those black-and-white bi-directed graphs.
There exists a homomorphism $f: G \to F$, if $F$ contains all vertices of
$G$, and all (black and white) edges of $G$, as  (black) edges of
$F$. White edges of $F$ are those which  start at the ends of
edges of $G$. So, these homomorphisms change white edges to a black ones
and add, as a white, new edges which come out from the former white
edges. Of course, the identity map from a black-and-white graph to itself is
considered to be an elementary homomorphism as well. Any composition
of the above defined elementary homomorphisms is also a homomorphism.
These homomorphisms play a role analogous to the 'complete graph
homomorphisms' in \cite{A}, however they are defined differently in
our case.

\begin{proposition}\label{prop1}\rm Every infinite, locally finite graph $E$ is a direct limit of a sequence of finite graphs and homomorphisms in the category $\cal E$. Any finite subgraph of $E$ can serve as a starting element of this sequence.
\end{proposition}

\begin{proof}
Take an arbitrary finite subgraph $E_0$ of $E$ (as a black-and-white subgraph constructed on vertices of $E_0$) and consider a sequence of non-identity homomorphisms $\phi_n: E_n \to E_{n+1}$ in $\cal E$. Due to the connectedness of the graph $E$, the union of edges of all elements of the sequence will coincide with the set of edges of $E$. Moreover, if an edge becomes black in the graph $E_n$ from the sequence, then it will be a black edge in any $E_N$, for $ N \geq n$. This means that $E$ is indeed a limit of a sequence $E_n, \phi_n$.
\end{proof}

After the category $\cal E$ is constructed, we define a functor from  $\cal E$ to abelian groups ${\cal AG}$.
We associate with a  black-and-white graph $E \in \cal E$ the group $K_0(\tilde {\cal O}_E)$ with generators corresponding to all (black and white) edges $x_i \in E$ and relations $x_i= \sum_{y_i \in E^1} \lambda_{i,j} y_j$, for any black edge $x_i \in E$. Here $y_j$ run over all edges, black and white of E, and $\lambda_{i,j}$ is equal to $1$ if there is a path connecting directly $x_i$ to $y_j$ (they are adjacent in the directed graph), except from the case when $y_i$ is the inverse of $x_i$. Otherwise $\lambda_{i,j}$ are equal to zero.

Let us describe how the functor ${\cal F}: {\cal E} \to {\cal AG}$ maps a black-and-white graph morphism $f: G \to F$ to the homomorphism of abelian  groups ${\cal F} (f): K_0({\cal O}_G) \to K_0({\cal O}_F)$. Since (due to the definition of graph morphisms) all edges of the graph $G$ (black and white) are also edges of the graph $F$ (black), and all relations of the group $K_0({\cal O}_G) $ are present in the group $K_0({\cal O}_F)$, the map $K_0({\cal O}_G) \to K_0({\cal O}_F)$ which sends generators of $K_0({\cal O}_G)$ to themselves, in $K_0({\cal O}_F)$, is a group homomorphism.

\begin{theorem}\label{contf}\rm

The above defined functor $\cal F$ from $\cal E$ to $\cal AG$ is continuous, i.e. it commutes with direct limits.
\end{theorem}

\begin{proof}
Elements of the direct limit of groups ${\cal G} = \underline{\lim} {\cal G}_n$ are sequences of elements of corresponding graphs ${\cal G}_n$, mapped to each other by corresponding morphisms (modulo the equivalence relation). In particular, generators of the limit group ${\cal G}$ are (classes of) sequences  $X=(x \to x \to x...)$ consisting of one generator $x$ of a particular group ${\cal G}_n$, with trivial maps. If we take into account the way how summation on these (classes of) sequences  is defined, we see that any relation on generators $X_1,...,X_n$ (any particular relation contains a finite number of generators) of the limit group ${\cal G}$
is actually present for generators $x_1,...,x_k$ of groups ${\cal G}_N, {\cal G}_{N+1},...$, starting from certain $N$.
In other words, any relation of the limit group $\cal G$ appears as a relation of some group ${\cal G}_N$, and stays the same in ${\cal G}_l$, for $l \geq N$. Our definition of group associated with the black-and-white subgraph was constructed in a way, which ensures that $K_0$ of an infinite graph has the set of relations, obtained as a union of relations in groups associated with finite subgraphs. This means that two sets of relations for  $\underline{\lim} K_0({\cal O}_{F_n})$ and for $K_0(\underline{\lim} {\cal O}_{F_n})$ are coincide.

\end{proof}

Combining this theorem with Proposition \ref{prop1} we have the following.

\begin{corollary}\label{cor1}\rm
Any $K_0$ group of a locally finite infinite graph $E$ is a direct limit of groups corresponding to finite subgraphs from the  category $\cal E$.
\end{corollary}

\section{$K_0$ calculations in the case of finite Betti number}

Now we turn to concrete calculations in the case when the Betti number of the graph is finite.

It is known  for finite graphs and row-finite graphs (see, for
example, \cite{CK} and \cite{PR}), that $K_0({\cal O}_E)= coker(Id-\Phi)$, where
$\Phi:{\mathbb Z}^{(E^1)} \to {\mathbb Z}^{(E^1)}$ is the homomorphism of countable direct sums of copies of $\mathbb Z$,
defined for the graph $E$ by the formula $(*)$.

First of all, we shall show that in any locally finite graph, we can perform any finite number of edge contractions, without changing $K_0$.

\begin{theorem}\label{blabla} Let the graph $E'$ be obtained from $E$ by contraction of one non-loop edge $x$ and its inverse $\bar x$. Then the groups $K_0({\cal O}_E)$ and $K_0({\cal O}_{E'})$ are isomorphic.
\end{theorem}

\begin{proof}
We will obtain the fact that the group $K_0$ is preserved under the edge contraction as a corollary of the following general lemma, which might be interesting in its own right.

\begin{lemma}\label{contr} Let $G$, $H$ be abelian groups and $T: G \oplus H \to G \oplus H$ a homomorphism, such that $Tx-x \in G$ for any $x \in H$. Let $P: G \oplus H \to G$ be a homomorphism such that $P|_G = id_G$ and $Px=x-Tx$, $x \in H$.

Then for $\tilde T: G \to G = P \circ T|_G$ the following is true:

$G \oplus H / T(G \oplus H ) \simeq G/ \tilde T(G)$

\end{lemma}

\begin{proof} (of lemma \ref{contr}).

Define a homomorphism $J: G/ \tilde T (G) \to G \oplus H / T(G \oplus H )$ as follows:

$$J(u+\tilde T (G)) = u+ T(G \oplus H ).$$

This map is well-defined, i.e. $\tilde T(G) \subseteq T(G \oplus H )$.
Let  $u \in G$ and $Tu=y+w, y \in H, w \in G $. Then $\tilde T u = P \circ Tu = P(y+w)=y-Tu+w=Tu-Ty \in T(G \oplus H)$.

The map $J$ is injective, i.e. for $u \in G, \,\, u \in T(G \oplus H)$ implies that $u \in \tilde T(G)$. Indeed, let $u=T(w+y) \in G$ for $w \in G, y \in H$. Since $u=Tw+Ty=Tw+y+Ty-y$, we can present $Tw$ as $Tw=(u+(y-Ty)) - y$, and $u+(y-Ty) \in G, y \in H$.
We showed above that for $w \in G, \,\, \tilde T w = T w = T h$, for $h$ being the $H$ component of $Tw$: $Tw=w'+h$, $w' \in G, h \in H$. Due to the above presentation of $Tw$ its  $H$ component is $-y$, so we have:
$ \tilde T w =Tw+Ty=T(w+y)$, hence indeed $T(w+y) \in \tilde T(G)$.

The map $J$ is surjective, i.e. $G + T(G \oplus H) = G \oplus H$.
Indeed, for $x \in H$, there exists $u \in G: u=Tx-x$. Then for $w+x \in G \oplus H$, $w+x=Tx-u+w$, where $w=u \in G$ and $Tx \in T(G \oplus H)$.

Thus $J$ is the required isomorphism.
\end{proof}

Now to prove Theorem \ref{blabla}, for an edge $x$ and its inverse $\bar x$ apply Lemma \ref{contr} for direct sum of copies of $\mathbb Z$: $G=\mathbb Z ^{(E^1 \setminus \{x, \bar x\})}$ and $H= \mathbb Z ^2$. As an operator $T$ (${\tilde T}$) we should take $T= Id - \Phi_E$ (${\tilde T} = Id - \Phi_{E'}$),  where $\Phi$ is defined by the formula (*).

\end{proof}

 We are now in a position to start the proof of the main theorem.

\begin{theorem}\label{t1} Let $E$ be an infinite, locally finite connected graph
with the finite first Betti number $\beta(E) $ and the valency
number $\gamma(E)$.

Then
$K_0({\cal O}_E)={\mathbb Z}^{\beta(E) + \gamma(E) }$.

\end{theorem}

\begin{proof}\hfill\break
{\it Type I}. If the valency number $\gamma(E)$ is finite, by a
finite number of steps we can reduce our graph to the rose with
$\beta(E)$ petals and $\gamma(E)$ outgoing simple infinite chains.
According to Theorem \ref{blabla}, $K_0$ will be preserved. In this
case it is easy to calculate directly the group $K_0({\cal
O}_E)=coker (Id-\Phi)$, generated by relations readable from the
graph. Indeed, let us denote variables corresponding to $\beta(E)=m$
petals (and their inverses) by $u_1,...,u_m, \
\bar{u}_1,...,\bar{u}_m$ and variables corresponding to
$\gamma(E)=n$ edges outgoing directly from the vertex of the rose
(and their inverses) by $x_1,...,x_n, \ \bar{x}_1,...,\bar{x}_n$.
Next, edges (and their inverses) in each chain will be
$x_i^{(1)},...,x_i^{(k)},...,\
\bar{x}_i^{(1)},...,\bar{x}_i^{(k)}..., \, i=\bar{1,n}$. Then $K_0$
will be the quotient of the free abelian group generated by the  set $\Omega
= \{ x_i^{(k)},\ \bar{x}_i^{(k)}, i=\bar{1,n}, u_j,\bar{u}_j,
j=\bar{1,m} \}$, subject to the relations defined by the formula
$K_0({\cal O}_E)=coker (Id-\Phi)$. For each edge $e \in E$ we will
have one relation. Note that the relation written for edges belonging to
chains will give $x_i^{(1)}=x_i^{(2)}=..., \, i=\bar{1,n}$. So after
that we actually have a finite number of relations for variables
$x_i,\bar{x}_i, \, i=\bar{1,n}, \, u_j,\bar{u}_j, \, j=\bar{1,m}$.
Namely,

$$ \sum_{j \ne k} (u_j+\bar{u}_j) + \sum_{l=1}^n x_l = 0 , \quad 1 \leq k \leq m,$$

$$ \bar {x}_l = \sum_{j=1}^m (u_j+\bar{u}_j) + \sum_{r \ne l} x_r, \quad 1 \leq l \leq n,$$

where the first group of relations corresponds to petals and the second to edges outgoing (incoming) from (to) the rose. It is a compete set of defining relations for $K_0$ on the set of generators $\Omega$.

For convenience, let us denote by $w_j=u_j+ \bar{u}_j$. Now write down the matrix of the above system of linear equations on variables $w_j, x_i,\bar{x}_i, \, j=\bar{1,m}, \, i=\bar{1,n}$.

$$
\left(
\begin{matrix}
%\begin{array}
%{ccccccccccccc
0&\dots& 1    &1&\dots&1   &0&\dots&0\\
 %1&\ddots& &1   &1&\dots&1   &0&\dots&0\\
 & \ddots&     & &\dots&    & &\dots & \\
1&\dots& 0    &1&\dots&1   &0&\dots&0\\
1&\dots& 1    &0&\dots&1   &-1&\dots&0\\
%1&\dots& &1    &1&\ddots&1  &0&\ddots&0\\
& \dots&     & &\ddots&    & &\ddots & \\
1&\dots& 1    &1&\dots&0   &0&\dots&-1
\end{matrix}
\right)
$$

By adding last $n $ columns to the first $m$ we can make zeros in the lower $n \times (m+n)$ block of the matrix. Then using the middle block of $n$ columns  we can transform the upper left $m \times (m+n)$ corner into

$$
\left(
%(a")
\begin{matrix}
1&0& &{\dots}&  & &0& 0&{\dots}&0  \\
0&1& &{\dots}&  & &0& 0&{\dots}&0    \\
& &{\ddots}& &  & & & & &   \\
0&1& &{\dots}&  & &1& 0&{\dots}&0    \\
\end{matrix}
\right)
$$

This shows that we have $m+n$ free variables: $u_1,...,u_m$,$\bar{x}_1,...,\bar{x}_n$.
So we see, that
$K_0({\cal O}_E)={\mathbb Z}^{\beta(E) + \gamma(E)} $ in this case.

{\it Type II}. The second case is when  the number $\gamma(E)$ is infinite. Here we can not write down a finite number of equations on the finite number of variables, which will define a group, but we can show what will be the system of free generators of the abelian group $K_0({\cal O}_E)$. The group $K_0({\cal O}_E)$ is defined by generators
corresponding to all edges of the graph, consisting of one rose with $\beta(E)$ petals and finite number
of outgoing infinite trees. The number of outgoing trees can not be infinite, because of the locally finiteness condition.

Let us consider generators corresponding to petals of the rose: $u_1,...,u_m,\\ \bar{u}_1,...,\bar{u}_m$ and edges coming out directly from the rose vertex: $x_1,...,x_n,\\ \bar x_1,..., \bar x_n$. We have the following equations on them:

$$ \sum_{k \ne j} (u_k+\bar{u}_k) + \sum_{l=1}^n x_l = 0 , \quad 1 \leq j \leq m$$

$$ \bar {x}_j = \sum_{k=1}^m (u_k+\bar{u}_k) + \sum_{l \ne r} x_l, \quad 1 \leq j \leq n$$

These are the same as above and
give us $n+m$ free variables: $u_1,...u_m,\\ x_1,...,x_n$. Then consider for any branching vertex, a piece of the tree of the shape

$$
%\left(
%(a")
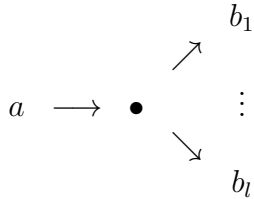
\begin{matrix}
& & & & {b_1} \\
& & &{\nearrow}& \\
a & {\longrightarrow} & {\bullet} & & {\vdots} \\
& & & {\searrow} & \\
& & & & {b_l}

\end{matrix}
%\right)
$$

Equations which we have to write for edges $a$, incoming for this vertex and $b_j, \, j=\bar{1,l}$,  outgoing from it, and leading to  an infinite path, form the following system.

$$ a = b_1 + ... + b_l$$

$$ \bar{b}_j = \bar{a} + {b}_1 + ... + \widehat{{b}_j}+...+{b}_l$$

So, on such a step we get $l-1$ new free variables,
corresponding to $l-1$ new infinite paths along the graph, we got in this vertex (which is equal to the valency of this branching vertex). If we sum up all new free variables, which we  got from  all vertices of outgoing trees, we arrive at $\gamma(E)$   additional variables. Note that again,  if an  outgoing chains  are finite, then
variables corresponding to their edges are just zero. So, we see that in this case also
$K_0({\cal O}_E)={\mathbb Z}^{\beta(E) + \gamma(E)}$, and here it is a direct sum of the countably infinite number of copies of $\mathbb Z$.
By this  the proof of the  theorem is completed.

\end{proof}

\section{The Whitehead group expressed via the first Betti number}

In the original paper due to Cuntz and Krieger \cite{CK} it was shown that $K_0$ and $K_1$ of  the Cuntz-Krieger $C^*$- algebra $ {\cal O}_A$, associated with any finite 0-1 matrix A
are, respectively, co-kernel and kernel of the map $ (Id - A^t): Z^n \to Z^n$.

This fact was later generalized in \cite{PR} to the
Cuntz-Krieger $C^*$- algebra $ {\cal O}_A$, associated in the same way with an infinite 0-1 matrix A, with a finite number of 'ones' in any row. The graph algebra we consider, as mentioned in the Introduction, is a
Cuntz-Krieger algebra associated with an infinite matrix, constructed from the graph by certain rules.
So this result is applicable here and making use of this  we will prove the following theorem.

\begin{theorem}\label{tk1} Let $E$ be an infinite, locally finite connected graph
with the finite first Betti number $\beta(E)$, and ${\cal O}_E$ is an associated (via the Bass-Hashimoto operator) $C^*$- algebra.
Then $K_1({\cal O}_E)={\mathbb Z}^{\beta(E)}$.

\end{theorem}

\begin{proof}

The proof will be divided into several steps.

\vspace{2mm}

First of all, we shall show that in any locally finite graph, we can perform any finite number of edge contractions, without changing $K_1$.

\begin{theorem}\label{blabla1} Let the graph $E'$ be obtained from $E$ by contraction of one non-loop edge $x$ and its inverse $\bar x$. Then the groups $K_1({\cal O}_E)$ and $K_1({\cal O}_{E'})$ are isomorphic.
\end{theorem}

\begin{proof}

Let us ensure the following fact of linear algebra.

\begin{lemma}\label{contr1} Let $G$, $H$ be abelian groups and $T: G \oplus H \to G \oplus H$ a homomorphism, such that $Tx-x \in G$ for any $x \in H$. Let $P: G \oplus H \to G$ be a homomorphism such that $P|_G = id_G$ and $Px=x-Tx$, $x \in H$.

Then for $\tilde T: G \to G = P \circ T|_G$ the following is true: ${\rm Ker} T \simeq {\rm Ker} \tilde T.$

\end{lemma}

\begin{proof} (of lemma \ref{contr1})

Take an element $u+y \in {\rm Ker} T$, with $u \in G$, $y \in H$. Let $Tu=w+x$, where $w \in G$, $x \in H$.
Then $\tilde T u = x - Tx + w = Tu - Tx$, so $Tu = \tilde T u + T x$. Now substituting that to $T(u+y)=0$, we have $0=T(u+y)=\tilde T u + T(x+y)$. From this we see first that $\tilde T u = - T (x+y)$.

Denote the $G$ and $H$ components
of an element $r \in G \oplus H$ by $r_G$ and $r_H$ respectively, so $r = r_G + r_H$, for $r_G \in G$ and $r_H \in H$.
Now comparing the $G$ and $H$ components of the left and right hand side of $\tilde T u =-T(x+y)=-(x+y)- g_{x+y}$, we have that $x+y=0$. Therefore $\tilde T u = - T(u+y)$ and $u+y \in {\rm Ker } T $  iff $u \in {\rm Ker } \tilde T $.
\end{proof}

The Proof of Lemma \ref{contr1} will follow as a corollary from this lemma if we put
$G=\mathbb Z ^{|E^1 \setminus \{x, \bar x\}|}$ and $H= \mathbb Z ^2$. As an operator $T$ (${\tilde T}$) we should take $T= Id - \Phi_E$ (${\tilde T} = Id - \Phi_{E'}$),  where $\Phi$ is defined by the formula (*).

\end{proof}

Theorem \ref{blabla1} from the first step, allows us to reduce the calculation of $K_1({\cal O}_E)$, where $E$ is a locally finite graph with the Betti number $\beta(E)$ to the $K_1$ for the graph $\Gamma$, which is a rose with $\beta(E)=\beta(\Gamma)$ petals
and a finite number of  trees, rooted in the vertex of the rose, with a finite number of branches outgoing from each vertex, due to the locally finiteness condition.

Now the second step in the calculation of the Whitehead group will be a calculation for the graph $\Gamma$.
We need to calculate $K_1({\cal O}_\Gamma) = {\rm Ker}\, T_{\Gamma}= {\rm Ker} (Id - \Phi_{\Gamma})$, where $\Phi_{\Gamma}$
is defined by formula (*).

For the graph $\Gamma$ we can present the set of all edges as a disjoint union of three sets:

$$ \Gamma^1 = \Gamma^{\uparrow} \sqcup \Gamma^{\downarrow} \sqcup R,$$

where  $\Gamma^{\uparrow}$ is the set of edges of the tree, directed towards the rose,
$\Gamma^{\downarrow}$ is the set of edges of the tree directed off the rose and $R$ is the set of petals of the rose.

Let $\xi \in {\rm Ker }\, T_{\Gamma}$, $\xi = \sum_{e \in \Gamma_{\xi}} m_e e$, where $\Gamma_{\xi}$ is a finite set of edges and $m_e \in \mathbb Z \setminus \{ 0 \}$.

Let us show first that the following is true.

\begin{lemma}
The set $\Gamma_{\xi}$ does not contain tree edges in the direction towards the rose: $\Gamma_{\xi} \cap \Gamma^{\uparrow} = \emptyset$.
\end{lemma}

\begin{proof} Consider the projection $\pi: \mathbb Z ^{(\Gamma)} \longrightarrow \mathbb Z ^{(\Gamma)} / \mathbb Z ^{(\Gamma')}$, where  $\Gamma' = \Gamma^{\downarrow} \sqcup R$, then denote by $T'$ the composition of our initial map $T$ with $\pi$:

$$\mathbb Z ^{(\Gamma)} \mathop{\longrightarrow}^T \mathbb Z ^{(\Gamma)} \mathop{ \longrightarrow}^{\pi} \mathbb Z ^{(\Gamma)} / \mathbb Z ^{(\Gamma')} \simeq \mathbb Z ^{(\Gamma^{\uparrow})}$$

Then $T' / \mathbb Z ^{(\Gamma')}=0$ and $T'e=e+f$, where $f$ consists of edges which are higher in the tree (=closer to the rose) than $e$.

Suppose $ \Gamma_{\xi} \cap \Gamma^{\uparrow} \neq \emptyset$.    Consider $g \in \Gamma_{\xi}$ farthest away from the rose, $\xi = mg + \tilde g, \, m \in \mathbb Z \setminus \{ 0 \}$.

Then

$$T' \xi = mg + f + T'(\tilde g).$$

Here $f$ consists of terms corresponding to the edges, closer to the rose than $g$.  $T'(\tilde g)$
consists of terms corresponding to the edges,  closer to the rose than $\tilde g$, which are in turn  closer then $g$. This means that the term $mg$ can not cancel, and $T'\xi \neq 0 $, hence $T \xi \neq 0$. We arrive at a contradiction.
\end{proof}

\begin{lemma}
The set $\Gamma_{\xi}$ does not contain tree edges in the direction off the rose: $\Gamma_{\xi} \cap \Gamma^{\downarrow} = \emptyset$.
\end{lemma}

\begin{proof}

Assume  $ \Gamma_{\xi} \cap \Gamma^{\downarrow} \neq \emptyset$, and take $h \in  \Gamma_{\xi} \cap \Gamma^{\downarrow}$, farthest 1 away from the rose. Then

$$ Th = h-(h_1+h_2+...),$$

where all $h_i$ are further away than $h$ from the rose.  Then if

$$ \xi=mh+\sum g_i,$$

 $m \in \mathbb Z \setminus \{ 0 \}$, then  $g_i \in \Gamma^{\downarrow}$, and since $g_i \in \Gamma^{\xi}$, they are closer than $h$.

  Therefore

$$T\xi = mh- m(h_1+h_2+...) + T (\sum g_i),$$

and the farthest from the rose edge, which could be contained in  $T(\sum g_i)$, is $h$. This means that the term $mh_1$
could not be cancelled, $T \xi \ne 0$ and we arrive at a contradiction.
\end{proof}

Now after the above two lemmas we are left with the only possibility, that

$$ \Gamma_\xi \subseteq  R = \{u_1,...,u_n, \bar u_1,..., \bar u_n \}.$$

Let $\xi= \sum_{j=1}^n (m_ju_j+n_j \bar u_j).$
We know that

$$T\bar u_j=Tu_j=-\sum_{k \neq j}(u_k+\bar u_k)+ \sum_{l=1}^s x_l=w-(u_j+\bar u_j),$$

 where for convenience we denote by $w=-\sum\limits_{k=1}^n (u_k+\bar u_k) + \sum\limits_{l=1}^s x_l.$
Then

$$ T \xi = (\sum_{j=1}^n (m_j+n_j))w - \sum_{j=1}^n (m_j+n_j)(u_j+\bar u_j)=0.$$

Since all $x_l$ appear in $w$ with coefficient 1 they can disappear only if $ \sum_{j=1}^n (m_j+n_j)=0.$
Hence

$$ T \xi = - \sum_{j=1}^n (m_j+n_j)(u_j+\bar u_j)=0.$$

Since each $u_j$ appears in one  term only, in order for the sum to be zero, we should have $m_j+n_j=0$ for all $j$.
This means that
$$\xi \in {\rm Ker} T \iff \xi=\sum_{j=1}^n m_j(u_j-\bar u_j),$$

 $m_j \in \mathbb Z$, thus
${\rm Ker} T = \mathbb Z^n$, where $n=\beta(\Gamma)$. This completes the Proof of Theorem \ref{tk1}.

\end{proof}

As a consequence of our result we have the following corollary.

\begin{corollary} An infinite analogue of the statement that $K_1$ is a torsion free part of $K_0$,
which holds for the case of finite graphs, is not true for infinite graphs.
\end{corollary}

\begin{proof}

It is well known, that in the case of finite graphs or matrices $K_1$ of a Cuntz-Krieger $C^*$-algebra associated
 with a matrix is a torsion free part of $K_0$ (see \cite{Ror} or \cite{MM}).

Let $E$ be, for example,  an infinite, locally finite graph with the
finite first Betti number $\beta(E)$ and infinite valency number
$\gamma (E)$. Then according to Theorem 3.3 $K_0({\cal O}_E) \simeq
\mathbb Z^{\infty}$ and $K_1({\cal O}_E) \simeq \mathbb
Z^{\beta(E)}$, according to Theorem 4.1. This shows that  $K_1({\cal
O}_E)$ is not isomorphic to the torsion part of $K_0({\cal O}_E)$.
\end{proof}

{\bf Acknowledgements}
It is my pleasure to express my gratitude  to the Max-Planck-Institut
f\"ur Mathematik in Bonn for hospitality and excellent research
atmosphere, where part of the work has been done. I am also thankful
to P.Ara for discussions at the early stage of the work on paper and
to the anonymous referee for many useful comments.
I would like to acknowledge the support from the grant FTE9038 of the Estonian Research Council.

\vskip1truecm

\scshape

\noindent  Natalia Iyudu \\

\noindent Max-Planck-Institut f\"ur Mathematik

\noindent 7 Vivatsgasse, 53111 Bonn

\noindent Germany

\noindent E-mail address: \qquad {\tt iyudu@mpim-bonn.mpg.de}\\

\noindent and
\noindent Queens's University Belfast

\noindent Department of Pure Mathematics

\noindent University road, Belfast, BT7 1NN, UK

\noindent E-mail address: \qquad {\tt n.iyudu@qub.ac.uk}

\vskip 5mm

\end{document}